\documentstyle[11pt,righttag]{amsart}

\begin{document}

\topmargin-0.1in
\textheight8.5in
\textwidth5.5in
 
\footskip35pt
\oddsidemargin.5in
\evensidemargin.5in

\newcommand{\V}{{\cal V}}      
\renewcommand{\O}{{\cal O}}
\newcommand{\LL}{\cal L}
\newcommand{\Ext}{\hbox{Ext}}
\newcommand{\pE}{\hbox{$^\pi\kern-2pt E$}}
\newcommand{\hQ}{\hbox{$\hat Q$}}
\newcommand{\phQ}{\hbox{$ '{\hat Q}$}}
\newcommand{\phd}{\hbox{$ '{\hat \delta}$}}

\newcommand{\lonto}{{\protect \longrightarrow\!\!\!\!\!\!\!\!\longrightarrow}}

\newcommand{\m}{{\frak m}}
\newcommand{\gl}{{\frak g}{\frak l}}
\newcommand{\ssl}{{\frak s}{\frak l}}
 
\renewcommand{\d}{\partial}

\newcommand{\ds}{\displaystyle}
\newcommand{\s}{\sigma}
\renewcommand{\l}{\lambda}
\renewcommand{\a}{\alpha}
\renewcommand{\b}{\beta}
\newcommand{\G}{\Gamma}
\newcommand{\g}{\gamma}

\newcommand{\C}{{\Bbb C}}
\newcommand{\N}{{\Bbb N}}
\newcommand{\Z}{{\Bbb Z}}
\newcommand{\ZZ}{{\Bbb Z}}
\newcommand{\K}{{\cal K}}
\newcommand{\ve}{{\varepsilon}}
\newcommand{\cupp}{{\star}}

\newcommand{\rowxy}{(x\ y)}
\newcommand{\colxy}{ \left({\begin{array}{c} x \\ y \end{array}}\right)}
\newcommand{\scolxy}{\left({\begin{smallmatrix} x \\ y
\end{smallmatrix}}\right)}

\renewcommand{\P}{{\Bbb P}}

\newcommand{\la}{\langle}
\newcommand{\ra}{\rangle}

\newtheorem{thm}{Theorem}[section]
\newtheorem{lemma}[thm]{Lemma}
\newtheorem{cor}[thm]{Corollary}
\newtheorem{prop}[thm]{Proposition}

\theoremstyle{definition}
\newtheorem{defn}[thm]{Definition}
\newtheorem{notn}[thm]{Notation}
\newtheorem{ex}[thm]{Example}
\newtheorem{rmk}[thm]{Remark}
\newtheorem{rmks}[thm]{Remarks}
\newtheorem{note}[thm]{Note}
\newtheorem{example}[thm]{Example}
\newtheorem{problem}[thm]{Problem}
\newtheorem{ques}[thm]{Question}
 
\numberwithin{equation}{section}

\renewcommand{\proofname}{\bf \rom{Proof}}
 
\newcommand{\onto}{{\protect \rightarrow\!\!\!\!\!\rightarrow}}
\newcommand{\donto}{\put(0,-2){$|$}\put(-1.3,-12){$\downarrow$}{\put(-1.3,-14.5) 

{$\downarrow$}}}

\newcounter{letter}
\renewcommand{\theletter}{\rom{(}\alph{letter}\rom{)}}

\newenvironment{lcase}{\begin{list}{~~~~\theletter} {\usecounter{letter}
\setlength{\labelwidth4ex}{\leftmargin6ex}}}{\end{list}}

\newcounter{rnum}
\renewcommand{\thernum}{\rom{(}\roman{rnum}\rom{)}}

\newenvironment{lnum}{\begin{list}{~~~~\thernum}{\usecounter{rnum}
\setlength{\labelwidth4ex}{\leftmargin6ex}}}{\end{list}}

\title{The Yoneda algebra of a $\K_2$ algebra need not be another $\K_2$ algebra}


\author[Cassidy, Phan and Shelton]{ }
\maketitle

\begin{center}

\vskip-.2in Thomas Cassidy$^\dagger$, Christopher Phan$^\ddagger$ and Brad Shelton$^\ddagger$\\
\bigskip

$^\dagger$Department of Mathematics\\ Bucknell University\\
Lewisburg, Pennsylvania  17837
\\ \ \\

   $^\ddagger$Department of Mathematics\\ University of Oregon\\
Eugene, Oregon  97403-1222
\end{center}

\setcounter{page}{1}

\thispagestyle{empty}

 \vspace{0.2in}

\begin{abstract}
\baselineskip15pt
The Yoneda  algebra of a Koszul algebra or a $D$-Koszul algebra is Koszul.  $\K_2$ algebras are a natural generalization of Koszul algebras, and one would hope that the Yoneda algebra of a $\K_2$ algebra would be another $\K_2$ algebra.  We show that this is not necessarily the case by constructing a monomial $\K_2$ algebra for which the corresponding Yoneda algebra is not $\K_2$.  
 \end{abstract}

\bigskip


\section{Introduction}
Let $A$ be a connected graded algebra  over a field $K$.  Correspondences between $A$ and its  bigraded Yoneda algebra $E(A)= \bigoplus_{n,m} Ext_A^{n,m}(K,K)$   have been studied in many contexts (e.g. \cite{G&M},  \cite{LPWZ2}, \cite{M-V} and \cite{Priddy}).   In particular there are very interesting classes of algebras where $E(A)$ inherits good properties from $A$.  Perhaps the most famous and intently studied of such classes of algebra is the class of Koszul algebras. 

An algebra  is  Koszul \cite{Priddy} if its Yoneda algebra is generated as an algebra by  cohomology degree one elements.  Koszul algebras will always have quadratic defining relations and given such an algebra, $A$, the Yoneda algebra is isomoprhic to the quadratic dual algebra $A^!$.  In particular, one has Koszul duality:  If $A$ is Koszul then $E(A)$ is Koszul and $E(E(A)) = A$.  

The following natural generalization of Koszul was introduced in \cite{CS-K2} and also investigated in \cite{Mauger} and \cite{Phan}.  We write $E^n(A)$ for $\bigoplus_p Ext_A^{n,m}(K,K)$.

\begin{defn} The graded algebra $A$ is said to be $\K_2$ if $E(A)$ is generated as an algebra by 
$E^1(A)$ and $E^2(A)$.
\end{defn}   

Koszul algebras are simply quadratic $\K_2$ algebras.  $\K_2$ algebras share many of the nice properties of Koszul algebras, including stability under tensor products, regular normal extensions and graded Ore extensions, (cf. \cite{CS-K2}).  Every graded complete intersection is a $\K_2$ algebra.   

Another important class of algebras is the class of $D$-Koszul algebras introduced by Berger in 
\cite{Berger}.  This is the class defined by: $Ext_A^{n,m}(K,K) = 0$ unless 
$m = \delta(n)$, where $\delta(2n) = nD$ and $\delta(2n+1) = nD+1$.    These algebras arise naturally in certain contexts and all such $D$-Koszul algebras are easily seen to be $\K_2$.   A remarkable theorem in \cite{GMMZ} states that if $A$ is $D$-Koszul algebra, then $E(A)$ is a $\K_2$ algebra, and furthermore, it is possible to regrade $E(A)$ in such a way that $E(A)$ becomes a Koszul algebra. In particular one gets a ``delayed'' duality:  $E(E(A)) = E(A)^!$ and $E(E(E(A))) = E(A)$. 

Based on the above theorem of \cite{GMMZ}, Koszul duality, and calculations of many other $\K_2$-examples, it seems reasonable to hope that the Yoneda algebra of any $\K_2$ algebra would also be 
$\K_2$, perhaps even Koszul.   Unfortunately, this is not always the case, and the purpose of this article is to exhibit an example of a $\K_2$ algebra for which the corresponding Yoneda algebra is not Koszul nor even $\K_2$.  Our example has 13  generators and 9 monomial defining relations. We believe that such a monomial algebra cannot be constructed with fewer generators and relations.
 
We wish to thank Jan-Erik Roos for pointing out an error in an earlier version of this paper.

\section{The algebras $A,  E(A)$ and $E(E(A))$}

Let $K$ be a field.  Let $\{m,n,p,q,r,s,t,u,v,w,x,y,z\}$ be a basis for a vector space $V$.  We define $A$ to be the $K$-algebra $T(V)/I$ where $I$ is the ideal generated by this list of monomial tensors:
$$R=\{mn^2p,\ n^2pqr,\  npqrs,\  pqrst,\  stu,\  tuvwx,\  uvwxy,\  vwxy^2,\  xy^2z\}.$$

\begin{thm}  The algebra $A$ is $\K_2$, but the algebra $E(A)$ is not $\K_2$.
\end{thm}

\begin{pf}
We use the algorithm given in section 5 of \cite{CS-K2} to prove that $A$ is $\K_2$.  From the set $R$ one can calculate that 
$S_1=\{m,n,p,q,r,s,t,u,v,w,x,y,z\}$,  $S_2=\{mn^2,\ n^2pq,\  npqr,\  pqrs,\  st,\  tuvw,\  uvwx,\  vwxy,\  xy^2\}$, $S_3=\{pqr, \ vw\}$,  $S_4=\{n^2\}$ and  $S_5= \emptyset$.  One easily verifies that for every $b\in S$ with minimal left annihilator $a$ we have either deg$(a)=1$ or $ab\in R$, and hence $A$ is $\K_2$. 
    
Let $B=E(A)$. In what follows we consider only the cohomology grading on $B$.  Following section 5 of \cite{CS-K2} we can construct a minimal projective resolution $P^\bullet$ of $_AK$ and see that the  Hilbert Series of the algebra $B$ is $1+13t^2+9t^2+8t^3+4t^4+3t^5+t^6$.   
 
It is possible (although laborious) to describe $B$ in terms of generators and relations and then construct a minimal resolution of $_BK$ and apply Theorem 4.4 of \cite{CS-K2} to show that   $B$ is not $\K_2$.  However $B$'s failure to be  $\K_2$ is apparent already in $Ext_B^3(K,K)$ and consequently there is a more efficient way for us to illustrate this.

Let $\bar m$ and $\bar z$ denote the basis elements in $B_1$ dual to $m$ and $z$ in $A_1$.  The vector space $B_2$ has a basis dual to the elements of the list of relations $R$.  We will use $\alpha, \beta$ and $\gamma$ to denote the dual basis elements  corresponding to the monomials $n^2pqr, stu$ and $vwxy^2$.  From the maps in the resolution $P^\bullet$ one can see that $\bar m \alpha$,  $\gamma\bar z$ and  $\beta \gamma$ are nonzero in $B$, while $\bar m \alpha\beta$ and $\beta\gamma\bar z$ are each zero. 

Recall that $Tor^B(K,K)$ can be calculated using the bar-complex \cite{PP} where ${\mathcal B}{ar}_i(K,B,K)=K\otimes_B\otimes B\otimes B_+\otimes \dots\otimes B_+\otimes K=B_+^{\otimes i}$.  Let $\zeta=\bar m\alpha\otimes\beta \gamma\otimes\bar z\in B_+^{\otimes 3}$. The differential on the bar-complex gives us $\partial(\zeta)=\bar m\alpha\beta \gamma\otimes\bar z -\bar m\alpha\otimes\beta \gamma\bar z=0$.  $\zeta$ is not in the image of $B_+^{\otimes 4}$ because $\partial (\bar m\otimes\alpha\otimes\beta \gamma\otimes\bar z)=  \bar m\alpha\otimes\beta \gamma\otimes\bar z-\bar m\otimes\alpha \beta \gamma\otimes\bar z$ while $\partial (\bar m \alpha\otimes\beta \otimes\gamma\otimes\bar z)= -\bar m \alpha \otimes \beta \gamma\otimes\bar z+\bar m \alpha \otimes \beta \otimes\gamma \bar z$.   Thus $\zeta$ represents a non-zero homology class in $Tor^B_3$.  

In contrast the element $\bar m\alpha\otimes\beta\gamma=\partial(-\bar m\alpha\otimes\beta\otimes\gamma)$ represents  zero in $Tor^B_2$ and $\bar m\alpha=\partial(\bar m\otimes\alpha)$   represents  zero in $Tor^B_1$.   
Therefore under the co-multiplication map $$\Delta:Tor_3^B(K,K)\to Tor_2^B(K,K)\otimes Tor_1^B(K,K)\oplus Tor_1^B(K,K)\otimes Tor_2^B(K,K)$$ we have $\Delta(\zeta)=0$.  This failure of $\Delta$ to be injective  is equivalent to the multiplication map $$E^2(B) \otimes E^1(B)\oplusÊÊE^1(B) \otimes  E^2(B) \to E^3(B)$$ not being surjective.  Hence $E(B)$ is not generated by $E^1(B)$ and $E^2(B)$, and so $B$ is not a $\K_2$ algebra.  \end{pf}

\bibliographystyle{amsplain}
\bibliography{bibliog}

\providecommand{\bysame}{\leavevmode\hbox to3em{\hrulefill}\thinspace}
\providecommand{\MR}{\relax\ifhmode\unskip\space\fi MR }
\providecommand{\MRhref}[2]{%
  \href{http://www.ams.org/mathscinet-getitem?mr=#1}{#2}
}
\providecommand{\href}[2]{#2}
\begin{thebibliography}{10}

\bibitem{Berger}
Roland Berger, \emph{Koszulity for nonquadratic algebras}, J. Algebra
  \textbf{239} (2001), no.~2, 705--734. \MR{MR1832913 (2002d:16034)}

\bibitem{CS-K2}
Thomas Cassidy and Brad Shelton, \emph{Generalizing the notion of {K}oszul
  algebra}, Math. Z. \textbf{260} (2008), no.~1, 93--114. \MR{MR2413345}

\bibitem{GMMZ}
E.~L. Green, E.~N. Marcos, R.~Mart{\'{\i}}nez-Villa, and Pu~Zhang,
  \emph{{$D$}-{K}oszul algebras}, J. Pure Appl. Algebra \textbf{193} (2004),
  no.~1-3, 141--162. \MR{MR2076383 (2005f:16044)}

\bibitem{G&M}
Edward~L. Green and Roberto Mart{\'{\i}}nez~Villa, \emph{Koszul and {Y}oneda
  algebras}, Representation theory of algebras (Cocoyoc, 1994), CMS Conf.
  Proc., vol.~18, Amer. Math. Soc., Providence, RI, 1996, pp.~247--297.
  \MR{MR1388055 (97c:16012)}

\bibitem{LPWZ2}
D.~M. Lu, J.~H. Palmieri, Q.~S. Wu, and J.~J. Zhang, \emph{A-infinity structure
  on ext-algebras}, pre-print, 2006.

\bibitem{M-V}
Roberto Mart{\'{\i}}nez-Villa, \emph{Skew group algebras and their {Y}oneda
  algebras}, Math. J. Okayama Univ. \textbf{43} (2001), 1--16. \MR{MR1913868
  (2003g:18013)}

\bibitem{Mauger}
Justin~M. Mauger, \emph{The cohomology of certain {H}opf algebras associated
  with {$p$}-groups}, Trans. Amer. Math. Soc. \textbf{356} (2004), no.~8,
  3301--3323 (electronic). \MR{MR2052951 (2005c:16014)}

\bibitem{Phan}
Christopher Phan, \emph{Generalized {K}oszul properties for augmented
  algebras}, preprint, 2008.

\bibitem{PP}
Alexander Polishchuk and Leonid Positselski, \emph{Quadratic algebras},
  University Lecture Series, vol.~37, American Mathematical Society,
  Providence, RI, 2005. \MR{MR2177131}

\bibitem{Priddy}
Stewart~B. Priddy, \emph{Koszul resolutions}, Trans. Amer. Math. Soc.
  \textbf{152} (1970), 39--60. \MR{MR0265437 (42 \#346)}

\end{thebibliography}

\end{document}